\documentclass[10pt]{amsart}
\usepackage{graphicx}
\usepackage{amsmath,amssymb}  
\usepackage{setspace}
\newtheorem{teo}{Theorem}
\newtheorem{defi}{Definition}
\newtheorem{prop}{Proposition}
\newtheorem{cor}{Corollary}
\newtheorem{lem}{Lemma}

\newcommand{\ZZ}{\mathbb{Z}}
\newcommand{\NN}{\mathbb{N}}
\newcommand{\RR}{\mathbb{R}}
\newcommand{\Fs}{\mathcal{F}^s}
\newcommand{\Fu}{\mathcal{F}^u}
\newcommand{\C}{\mathcal{C}}
\newcommand{\diff}{\text{Diff }}

\renewcommand{\SS}{\mathbf{S}}
\newcommand{\Cyl}[1]{\mathbf{C}_{#1}}
\renewcommand{\sl}{\text{sl}}
\newcommand{\eps}{\epsilon}

\title{$C^k$-Robust transitivity for surfaces with boundary}
\author{Aubin Arroyo}
\author{Enrique R. Pujals}

\thanks{MSC-2000 Classification:  37E30; 37G25; 37D50}
 
\address{A. Arroyo. Instituto de Matem\'{a}ticas, Unidad Cuernavaca.
 Universidad Nacional Aut\'{o}noma de M\'{e}xico,
A.P. 273-3 Admon. 3, Cuernavaca, Morelos, 62251 M\'{e}xico}
\email{aubin@matcuer.unam.mx}
\address{E. R. Pujals. Instituto Nacional de Matem\'atica Pura e Aplicada.
 Estrada Dona Castorina 110, CEP  22460-320.  Rio de Janeiro, Brasil }
\email{enrique@impa.br}
 
\begin{document}
  
\begin{abstract}
We prove that $C^1$-robustly transitive diffeomorphisms on surfaces with boundary do not exist, and we exhibit a class of diffeomorphisms of surfaces with boundary which are $C^k-$robustly transitive, with $k\geq 2$.
This class of diffeomorphisms are examples where a version of Palis' conjecture on surfaces with boundary, about homoclinic tangencies and uniform hyperbolicity, does not hold in the $C^2-$topology.
This follows showing that  blow-up of pseudo-Anosov diffeomorphisms on surfaces without boundary, become $C^2-$robustly topologically
mixing diffeomorphisms on a surfaces with boundary. 
\end{abstract}
\maketitle

\section{Introduction}
A diffeomorphism on a compact manifold $M$ is said $C^k-$robustly transitive if every diffeomorphism in a $C^k-$neighborhood of it has a dense orbit in $M$, for $k\geq 1$.  
The description of these systems is an important challenge: on one hand, being robust, they can not be ignored in any global 
picture of dynamical systems; on the other hand, they often exhibit a chaotic dynamical behavior.
Typical models showing robust properties for boundaryless compact manifolds
are the well known Anosov maps, which are uniformly hyperbolic in the whole manifold. 
In the case of surfaces without boundary, the complete picture is described in the $C^1-$topology, by the result of Ma\~n\'e, in \cite{Ma}, which states that
any $C^1-$robust transitive diffeomorphisms is Anosov. 

Surprising new results appear when one considers manifolds with boundary. 
It was shown by MacKay in \cite{Mac}, that a steady mixing (even Bernoulli) smooth volume-preserving vector field in a bounded container in $\RR^3$ with smooth no-slip boundary exists in a robust way: any (conservative) $C^3-$perturbation keeps the same properties. In particular, the flows introduced in \cite{Mac} are $C^3-$robustly transitive.

Let us call by a \emph{surface with boundary} a compact orientable riemaniann surface with boundary, and denote it by $(S,\partial S)$.
The boundary is $\partial S$ and consists in a finite union of disjoint closed smooth curves.
The genus of $(S,\partial S)$ is the genus of the surface obtained by collapsing each boundary component of $S$ into a point.
Let us denote by $\diff^{k}(S,\partial S)$ the set of diffeomorphisms of class $C^k$ defined on a given surface with boundary $S$.
Main theorem in this paper is the following.

\begin{teo} \label{cedos}
For any surface with boundary $(S,\partial S)$ with genus larger than zero, the set of $C^2$-robustly transitive diffeomorphisms 
is not empty. If the genus of $(S,\partial S)$ is zero and $\partial S$ has at least four connected components, the same result is true. 
\end{teo}

It is remarkable that we cannot obtain a similar result for the $C^1-$topology. In fact, here we prove the following theorem.

\begin{teo} \label{ceuno}
Given a surface with boundary $(S,\partial S)$ the set of $C^1-$robustly transitive diffeomorphisms in $\diff^1(S,\partial S)$ is empty.
\end{teo}

Examples on Theorem \ref{cedos} are constructed through a blow-up procedure of some fixed points of a Thurston's pseudo-Anosov map defined on a particular boundaryless surface (see Theorem \ref{p-a-thm}). 
On the other hand, Thurston's pseudo-Anosov maps exist with any prescribed set of singularities, except for the case of the sphere, where some topological obstruction appear.
For the case of the sphere with boundary, if the number of connected components of $\partial S$ is larger than four, this obstruction vanishes.

\subsection{Blow-up of pseudo-Anosov maps}
Let $S_{0}$ be a surface without boundary. Following Thurston's definition (see  \cite{Thurston}, \cite{FathiLaudembach}, \cite{GerberKatok}), a homeomorphism $f:S_{0}\to S_{0}$ is a pseudo-Anosov map if there is a finite set $\text{Sing}(f) \subset S_{0}$ such that $f$ is $C^\infty$ in $S_{0} \smallsetminus \text{Sing}(f)$, there are two two measured foliations $(\Fs, \mu_s)$ and $(\Fu, \mu_u)$ which are $f$-invariant, and such foliations are transversal each other in $S_{0} \smallsetminus \text{Sing}(f)$. Moreover, $f$ preserve a natural absolutely continuous measure $\mu = \mu_s \times \mu_u$, whose density is $C^\infty$ except at singular points, and $f$ is Bernoulli with respect to this measure. 
Also, there is a number $ \lambda >1$  such that:
\begin{equation*} 
f_*(\mu_s) = (1/ \lambda) \mu_s \text{ and } f_*(\mu_u) =   \lambda  \mu_u,
\end{equation*}
where $f_{*}(\nu )$ is the push-forward of the measure $\nu$ restricted to the corresponding foliation. The number $\lambda$ is called the dilatation of $f$ and $\log(\lambda)$ is precisely its topological entropy. 
Any point in $\text{Sing}(f)$ is a $p$-prong singularity for some $p \geq 1$, and is contained in $\text{Fix}(f)$, the set of fixed points of $f$.
These singular points are not hyperbolic, except for $p=2$. However, any periodic point of $f$ which is not a singular point is hyperbolic. 
 
A blow-up of surface at a point $x$ is, roughly speaking, the construction of a new surface from the given one, replacing only one chart at $x$ by another at some circle at the boundary. This replacement is done using a polar change of coordinates, and can be performed successively on a finite number of different points. There is a correspondence between several objects (\emph{e.g.} vector fields, maps) of a surface and its blow-up. This is precised in Lemma \ref{blowupconstruction} of Section \ref{partedos}.
 
\begin{defi}
Let $S_{0}$ be a surface without boundary, let $f_{0} : S_{0} \to S_{0}$ be a pseudo-Anosov map, and let $B$ be a finite subset of $S_{0}$ that $\text{Sing}(f_{0})\subset B \subset \text{Fix}(f_{0})$. A blow-up of a pseudo-Anosov map is the map $f \in \diff^{\infty}(S, \partial S)$, obtained after blowing-up every point in $B$.
\end{defi}

In Section \ref{partedos} we shall prove the following Theorem:
\begin{teo}\label{p-a-thm}
Any blow-up of a pseudo-Anosov map is $C^2-$robustly transitive. 
\end{teo}

In the classical work \cite{FathiLaudembach}, Fathi, Laudembach and Po\'enaru prove that there exist transitive pseudo-Anosov maps with only one singular point, on any surface without boundary of genus $g>1$. On the torus, any Anosov map is pseudo-Anosov.  
In \cite{MasurSmillie}, Masur and Smillie established that the existence of pseudo-Anosov maps in a boundaryless surface $S_{0}$ is strongly related to the existence of quadratic differentials on the surface. In fact, they guarantee their existence on every orientable surface without boundary with any prescribed set of singularities, if the number of prongs, at each singular point, satisfy certain relation with the genus of the surface.  For our purposes, we can understate such theorem in the following way: 

\begin{teo}[Masur-Smillie] \label{ms}
If a finite set $\, \{(x_{i}, p_{i} ) \in S_{0}\times \NN \, | \, 1 \leq i \leq k\}$ satisfy that:
\begin{equation} \label{mscondition}
\sum_{ i=1}^k(p_{i}-2) = 4(g-1),
\end{equation}
then there is a pseudo-Anosov homeomorphism $f :S_{0}\to S_{0}$ with orientable stable foliation, such that $ \text{Sing}(f)=\{x_{1}, \ldots, x_{k} \}$ and, for $i=1,\ldots,k$, the number of prongs at $x_{i}$ is $p_{i}$.
\end{teo}
Notice that $2$-prong singularities correspond to hyperbolic fixed points, and if $p_{i}=2$ for some indexes, they does not affect the left side of the equation (\ref{mscondition}). 
On the other hand, given a surface $S$ of genus $g\geq 1$ and $m\geq 1$ connected components of $\partial S$, Theorem \ref{ms} guarantees the existence of a pseudo-Anosov map on a  surface without boundary of genus $g$, with enough fixed points to rebuild the boundary after several blow-ups. 
Therefore, Theorem \ref{cedos} follows directly from Theorem \ref{p-a-thm}.

The topological obstruction for the case of the sphere is that any pseudo-Anosov map on the sphere exhibits at least four 1-prong singularities. 
Examples of these maps can be obtained from a linear Anosov map of the torus, identifying the sphere with the quotient space of the torus by the involution.
Observe also that, for $g=0$, the simplest configuration to obtain (\ref{mscondition}) is with four points with $p_{i}=1$. If we include some other points with $p \geq 3$, it has to be compensated by several 1-prongs.
 
There is a general definition of a pseudo-Anosov map on surfaces with boundary. In fact, it is well established that this definition
has to allow that the diffeomorphism is Morse-Smale restricted to the boundary, and not the identity map, as first stated in \cite{FathiLaudembach} (see \cite{GerberKatok}). 
Recall that a map of a one dimensional compact manifold (may not be connected) is Morse-Smale if the Limit set consist only of a finite number of hyperbolic periodic points.
However, for our purposes, it is not clear that we can find adequate smooth charts of a general pseudo-Anosov map around the boundaries.

Given a hyperbolic periodic point $q$, denote the homoclinic class of $q$ by:
$$H_{q}  =  \{ p \in \text{Per}(f) \, | \, W^s(p) \cap W^{u}(q) \neq \emptyset \text{ and }  W^u(p) \cap W^{s}(q) \neq \emptyset\},$$
where  $W^s(p)$ and  $W^u(p)$ are the stable and unstable manifolds of a hyperbolic periodic point. 
If a diffeomorphism $f$ on a compact manifold has a periodic point whose homoclinic class is dense in the whole space, then it is topologically mixing;
that is, for any pair of not empty open sets $A$ and $B$ there exist an integer $m\geq 0$ such that for any $n\geq m$ holds that $f^n(A)\cap B\neq \emptyset$.
If this happens for any $g$ in a $C^k-$neighborhood of $f$, it is said that $f$ is $C^k-$robustly topologically mixing.
Related to that, we prove the following. 
\begin{teo} \label{cedos-hc}
If $g \in \diff^2(S, \partial S)$ is $C^2-$close to a blow-up of a pseudo-Anosov map, then there is a hyperbolic periodic point of $g$ whose homoclinic class is dense in $S$. 
\end{teo}

All these results suggest that these systems may be $C^2-$structurally stable, however, some work has to be done to prove or disprove it.
Finally, it is natural to wonder if a blow-up of a pseudo-Anosov map
have a unique physical measure. 
It is worth to mention that hyperbolicity, by itself, is not sufficient to understand global robust properties in the context of surfaces with boundary. 
In fact, there are no Anosov maps on surfaces with non empty boundary.   
Besides that, it is not known if there are examples of $C^2-$robust systems which are not robust in the $C^1-$topology, on surfaces without boundary.

\subsection{Version of Palis' conjecture for surfaces with boundary}
From another point of view, the characterization of robust dynamics provide a conceptual scheme
that helps to describe the dynamical behavior of the \emph{generic} dynamical system. 
From this perspective, in the early 80's,  Palis conjectured for boundaryless surfaces the following (see, for instance, \cite{palis}):

\vspace{8pt} \noindent \textbf{Conjecture} (Palis).
\emph{Given $k \geq 1$, every $C^k$-diffeomorphism of a compact surface without boundary can be approximated, in the $C^k-$topology, by one which is hyperbolic or by one exhibiting a homoclinic tangency.} 
\vspace{8pt} 

Recall that  a homoclinic tangency between two hyperbolic periodic points $p$ and $q$ is a point $z \in W^s(p) \cap W^{u}(q)$, and where this intersection is not transversal. 
In \cite{PS2}, this conjecture is proved to be true  in the $C^1-$topology, for both cases: with or without boundary. 
In this paper we obtain that this conjecture is false when is formulated for compact surfaces with not empty boundary, in the $C^2-$topology. This is a  consequence of Theorem \ref{cedos}, and follows immediately from the next theorem:

\begin{teo} \label{nonhyp}
If $g \in \diff^2(S, \partial S)$ is $C^2-$close to a blow-up of a pseudo-Anosov map, then $g$ is not hyperbolic and it do not exhibit homoclinic tangencies between periodic points.
\end{teo}
 
In this direction, a blow-up of a pseudo-Anosov map can be suspended into a non-singular flows in certain three dimensional manifolds with boundary.
These flows are counter-examples for the natural generalization of Palis' conjecture to non-singular flows on $3-$manifolds with boundary, in the $C^2-$topology. Even though, it is proved in \cite{ARH} a positive statement of this conjecture, in the $C^1-$topology, which includes flows with singularities. 
\begin{cor}
There are examples of smooth flows without singularities, defined on three dimensional manifolds with boundary, which are not approximated, in the $C^2-$topology, by neither hyperbolic ones nor ones exhibiting a homoclinic tangency.
\end{cor}

\subsection{About the proofs of Theorem \ref{ceuno} and \ref{p-a-thm}}
Topological rigidity of the boundary is one of the key facts that allows us to prove robust transitivity for blow-up of pseudo-Anosov maps. In fact, we found that these maps have a non-uniformly hyperbolic structure that persists under $C^2-$perturbations. This property is related to the notion of \emph{hyperbolic cone structure} for boundary maps defined in Section \ref{boundarymaps}, and provides also good estimates on the growth of lengths of curves inside certain cone-fields, when they are close to the boundary. 
In order to prove that any small $C^2-$perturbation of these maps are transitive, we mix this last property with some robust properties, in the $C^1-$topology, which are valid only outside of a neighborhood of the boundary. 
As we shall show, these maps inherit a stable and unstable transversal foliations at any point of the surface. However, the angle between them is not bounded away from zero. 
Nevertheless, this structure is not preserved for $C^1-$perturbations.  
In fact, a small $C^1$-perturbation can create a tangency and therefore either a periodic sink or repeller, see \cite{BonattiDiazPujals}.
This argument is the key to prove Theorem \ref{rtsets}. 
In Section \ref{parteuno} we give a proof of a more general statement: any non-trivial robustly invariant sets can not intersect the boundary of the surface, otherwise it is possible to create a sink or a source by an arbitrarily small $C^1-$perturbation of the original map.  

\subsection{Acknowledgment} We thank Adolfo Guillot for helpful conversations in the preparation of this paper. First author was partially supported by CONACyT grant 58354 and PAPIIT-UNAM grant IN102307.


\section{$C^1$-Robustly transitive sets} \label{parteuno}
The proof of Theorem \ref{ceuno} relies on a more general concept of robust transitivity focused not in the whole space but on maximal invariant sets of a diffeomorphism. Let $(S, \partial S)$ be a surface,  $f \in \text{Diff }^1(S, \partial S)$, and let $U$ be an open subset of $S$.
Assume that $\Lambda_f(U)= \cap_{n\in \ZZ}f^n(U) \neq \emptyset$ is a compact transitive invariant set of $f$. 
We say that the set $\Lambda_f(U)$ is a $C^1$\emph{-robustly transitive set} if there is an $C^1$-open neighborhood $\mathcal{N}$ of $f$, such that for any $g \in \mathcal{N}$ the set $\Lambda_g(U)= \bigcap_{n\in \ZZ}g^n(U)$ is a transitive compact invariant set of $g$.

\begin{teo}  \label{rtsets}
Let $\Lambda$ be a $C^1-$robust transitive set that $\Lambda  \cap \partial S \neq \emptyset$, then $\Lambda$ is a hyperbolic fixed point in the boundary.
\end{teo}

\begin{proof}
Let us assume that $\Lambda$ is not a single point and that $\Lambda \cap \partial S \neq \emptyset$. 
Let $B$ be connected component of the boundary that intersect $\Lambda$.
Observe that generically $f|_{B}$ is Morse Smale systems, so, from the fact that we assume that $\Lambda$ is not reduced to a point, it follows that $\Lambda \cap (S\smallsetminus \partial S) \neq \emptyset$.
Moreover, there are two saddle-type fixed point $p$, $q$ in $S$ accumulated by $\Lambda$. 
We assume that $W^u(p)\subset B$, $W^s(q)\subset B$ and that there is a saddle connection between these two points inside $B$.
Let us denote by $\gamma$ the arc inside $B$ given by $\gamma := W^u(p)\cap W^s(q)$ and observe that the extremal points of this arc are given by $p$ and $q$. 
Notice that also holds that $W^s(p)\subset (S\smallsetminus \partial S)$ and  $W^u(q)\subset (S\smallsetminus \partial S)$. 
From the fact that we are dealing in the $C^1-$topology, we can assume that $\Lambda$ is a homoclinic class and there exists a periodic point $\hat p$ such that $W^u(\hat q) \cap W^s(p) \neq  \emptyset$ and 
$W^s(\hat q) \cap W^u(p)\neq  \emptyset$. 
Therefore, there are points in the homoclinic class of $\hat q$ accumulating on $\gamma$ and so, it follows that the angle between the stable manifold and the unstable manifold of $\hat q$ at those points goes to zero. 
From that, using an argument as the one given in \cite{BonattiDiazPujals}, it follows that by a $C^1$ perturbation it is created either a sink or a repeller. 
Hence, $\Lambda$ is a hyperbolic fixed point in the boundary.
\end{proof}

\begin{teo}
Let $\Lambda$ be a non trivial $C^1-$robust transitive set. Then, it follows that $\Lambda$ is hyperbolic and $\Lambda \cap \partial S = \emptyset$.
\end{teo}

\begin{proof} Let $\Lambda$ be non trivial $C^1-$robust transitive set. Theorem \ref{rtsets} implies that $\Lambda \cap \partial S = \emptyset$. Therefore, we are reduced to the case of boundaryless surfaces, where it is well known that $C^1$-robustly-transitive sets are  hyperbolic.
\end{proof}


\section{$C^2$-robust properties of boundary maps} \label{partedos}
In this section we give a proof of Theorem \ref{p-a-thm} and therefore, we conclude Theorem \ref{cedos}. 
To prove this Theorem we need to define the blow-up of a a pseudo-Anosov map, introduce the notion of boundary maps and study carefully the local model of this maps in a neighborhood of a component of the boundary.
Denote the circle by $\SS^1$.
The following lemma is the key to construct a blow-up of a pseudo-Anosov map. 
\begin{lem} \label{blowupconstruction}
Let $S_{0}$ be a smooth surface, and $\sigma \in S_{0}$. Let $f_{0}:S_{0}\to S_{0}$ be a pseudo-Anosov map on $S_{0}$ that $\sigma \in \text{Sing}(f_{0}) \cup \text{Fix}(f_{0})$. There is a surface $S$ with boundary and a map $f \in \text{Diff}^\infty(S, \partial S)$ that is $C^{\infty}-$conjugated to $f_{0}$ in $S_{0}\smallsetminus \{\sigma\}$, and there is a chart around $\partial S$ where $f$ has the following expression:
\begin{equation} \label{themap2}
\left( \frac{2}{p}  \arctan( \lambda^2 \tan(  px/2   ) ) , y \sqrt{ \frac{\cos^2(px/2)}{\lambda^{2} }  + \lambda^2 \sin^2(px/2) } \right).
\end{equation}
\end{lem}

\begin{proof}
Let $S_{0}$ be a smooth surface and let $\sigma \in S_{0}$. Consider an atlas of $S_{0}$, that is, a collection of charts $\{ ( U \subset S_{0}, \phi:U \to \RR^2 ) \}$,
such that the change of coordinates, say $\psi \circ \phi^{-1}$, are $C^{\infty}$. 
Take a chart $(U_{\sigma}, \phi_{\sigma})$ at $\sigma$. 
The blow-up over the point $\sigma$ is the topological space 
$ S  = (S_{0} \smallsetminus \{\sigma \}) \cup \SS^1,$
provided with a new atlas formed by all charts of $S_{0}$ on points of $S_{0} \smallsetminus \{\sigma \}$,
and, an additional chart $(U^*, \phi^*)$, obtained by the expression of $\phi_{\sigma}$ in polar coordinates $(x, y )$, where $x$ is the angle and $y$ is the lenght. Notice that $S$ is a smooth surface with boundary.

Now let $f_{0}$ be a pseudo-Anosov map on $S_{0}$ and assume that  
$\sigma \in \text{Sing}(f_{0}) \cup \text{Fix}(f_{0})$, with $p \geq 1$ prongs. Recall that if $p=2$, the point $\sigma$ is hyperbolic.
The unstable prongs landing at $\sigma$, in a small neighborhood, determine $p$ radial sectors, $\Gamma_{1}, \ldots, \Gamma_{p}$, which we can assume are fixed by $f$. 
Denote by $B_{r}(0,0) \subset \RR^2$, the ball of radius $r>0$ around the origin. 
There is a $C^{\infty}-$chart $\phi:U\to B_{r}(0,0) \subset \RR^2$, for some $r>0$ and for some neighborhood $U$ of $\sigma$, (see \cite{GerberKatok}, for instance).
In this chart the map $f|_{\Gamma_{i}\cap U}$, for $i \in \{1, \ldots p\}$ can be written, in polar coordinates, as equation (\ref{themap2}). 
Notice that this map extends to a $C^\infty$ map to the line $[y=0]$.
It is easy to see that  it induces a diffeomorphism  $f \in \text{Diff}^\infty(S,\partial S)$. 
\end{proof}

The topology of the surface $S$ is independent of the choice of $x$. In particular, the genus of $S$ is the same of $S_{0}$.
Moreover, the surface $S$ obtained is equivalent to the one obtained as the connected 
sum of the original surface and a real projective space of dimension 2, and then cutting along the unique not null-homotopic curve of the projective space. This curve corresponds to the new connected component of the boundary of $S$.
In this way, any surface $S$ of genus $g$ and with $m\geq 1$ components of the boundary can be obtained from boundaryless surface $S_{0}$ of genus $g$ after a finite number of blow-ups on different points $\{x_{1}, \ldots , x_{m}\} \subset S_{0}$.

\subsection{Cone-fields}
Let $E$ be a two-dimensional real vector space. 
A \emph{half-cone} $C$ is a proper subset $\{0\} \neq C\subset E$ such that: $C+C \subset C$; $t \cdot C \subset C$, for any $t\geq 0$; and $C\cap (-C) = \{0\}$. A (complete) cone is $C \cup (-C)$, for a given half-cone $C$.  
A standard way to measure the angle between two vectors $v, w \in E$ is with the interior product: 
\[ \angle(v,w) = \arccos \left( \frac{<v,w>}{||v|| \, ||w||} \right).\] 
Let $\{\hat e_{1}, \hat e_{2} \}$ be the canonical basis of $E$. 
Given a non-zero vector $v = (v_{1},v_{2}) \in E$, the angle between $v$ and $\hat e_{1}$ can be computed in terms of the slope of $v$, that is:
$ \sl(v) = v_{2}/v_{1}$. In fact, $\angle(v, \hat e_{1}) = \arctan (\sl(v))$. In the same way, $\sl^{\perp}(v) = v_{1}/v_{2}$ measures the angle between $v$ and $\hat e_{2}$.

\begin{figure}
\includegraphics*[width=5in]{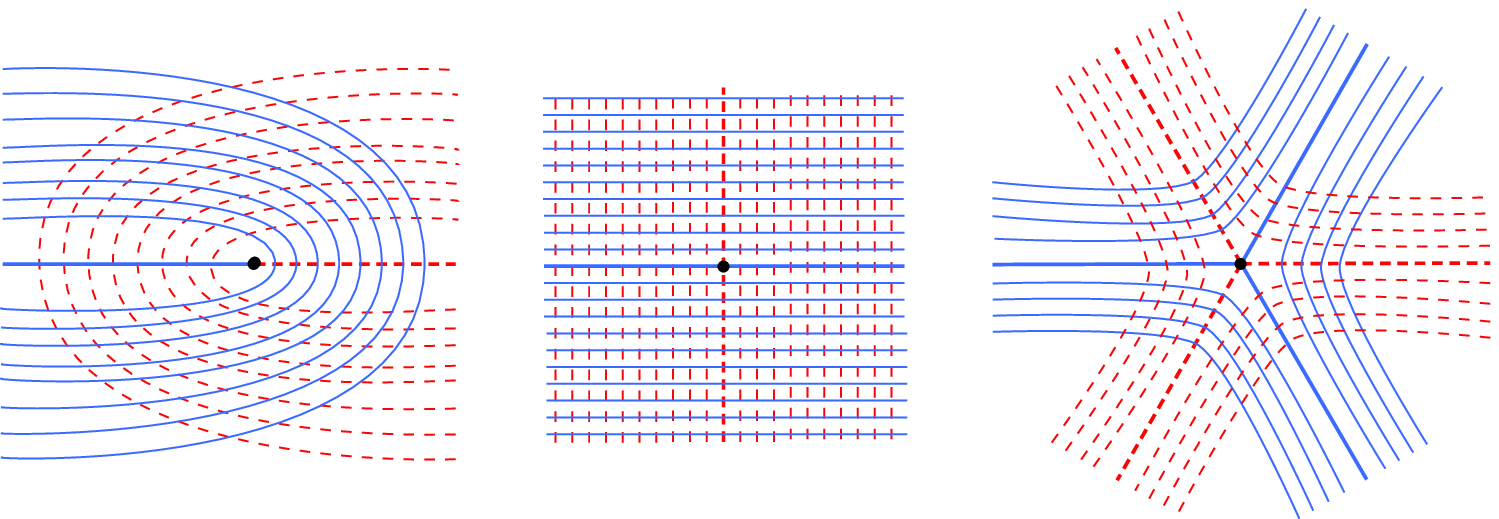}
\caption{A $p$-prong singularity of $\Fu$ and $\Fs$ with $p=1,2$ and $3$.}
\label{fig:prong}
\end{figure}

Let $f$ be a pseudo-Anosov map on a surface without boundary $S_{0}$ of genus $g \geq 0$. 
By definition, stable and unstable foliations of $f$ have a finite number of common $p$-pronged singularities with $p\geq 1$. Let $x$ be a singular point or a hyperbolic fixed point of $f$. There is an integer number $p_{x}\geq 1$ of leaves of $\Fu$ landing on $x$ and the same number of leaves of $\Fs$, as depicted in figure \ref{fig:prong}. 

In \cite{GerberKatok} is proved that
there are two invariant cone-fields in $S_{0} \smallsetminus \text{Sing}(f)$, along $\Fu$ and $\Fs$ respectively, defined by:
$$K^{u,s}(\alpha; x) = \{v \in T_{x}S \, | \, \angle(v, T_{x}\mathcal{F}^{u,s}) < \alpha \},$$ 
for some $0 < \alpha < 1$, and these cones satisfies that: 
\begin{equation*} 
D f  K^{u}(  x) \subset K^{u }(  f(x)) \text{ and } Df^{-1} K^{s}(  x) \subset K^{u }(  f^{-1}(x)).
\end{equation*}
If $W_{0}$ is a neighborhood of $\text{Sing}(f)$, there is $\tilde \lambda > 1$ such that for any point $x \in S_{0} \smallsetminus W_{0}$ we have that:
$|D f_{x}(v)| > \tilde \lambda |v|$, for any $v \in K^{u}(\alpha, x)$ and $|D f^{-1}_{x}(v)| > \tilde \lambda |v|$, for any $v \in K^{s}(\alpha, x)$.
Moreover, there is $\zeta > 0$ such that for any  $x \in S_{0} \smallsetminus W_{0}$ we have that:
$$\angle (K^{s}( x), K^{s}( x )) > \zeta.$$
If the angle between the stable and unstable cones are bounded away from zero, 
there is a $C^1$-neighborhood $\mathcal{N}_{1}$ of $f$ such that the same cones $K^{u,s}( \, \cdot \,)$ are invariant under $D  g$, and expanded and contracted by the same constant $\tilde \lambda>1$, for any $g \in \mathcal{N}_{1}$.

\begin{lem} \label{unifhyp}
Any maximal invariant set of $g \in \mathcal{N}_{1}$ contained in $S \smallsetminus W_{0}$ is uniformly hyperbolic. Moreover, it is contained in the homoclinic class of certain hyperbolic periodic point.
\end{lem}

\begin{proof}
The bound on the angle follows from the fact that $S \smallsetminus W_{0}$ is a compact set and the fact that the foliations $\Fu$ and $\Fs$ are transversal 
in $S \smallsetminus W_{0}$. Uniform hyperbolicity for such a set is obtained by standard arguments. 
\end{proof}

\subsection{Local model around the boundary} \label{boundarymaps}
In this section we shall focus on a family of maps of the cylinder $\Cyl{r} = \SS^1 \times [0,r) \subset \RR^2$, for some $r>0$,
that corresponds to the blow-up of a pseudo-Anosov map at a singular point,
in a neighborhood of a component of the boundary. 
In particular, we are interested in the existence of input and output sections which intersects any orbit that accumulates on such boundary component, and in the existence of two invariant cone-fields for $C^2-$perturbations of the transition map between these sections. 

Let $R^+$ and $R^-$ be two open and disjoint open sets of $\Cyl{r}$ and let $T:R^+ \to R^-$ be a piecewise smooth map.
We say that $T$ has a \emph{hyperbolic cone structure} if there are numbers $0< \alpha_0 < \alpha<1$, $\sigma >1$,  and there are two continuous cone-fields in $R^+ \cup R^-$, say $C^{u}( \alpha;x )$ and $C^{s}( \alpha;x )$, such that:
\begin{itemize}
\item  $C^{u}( \alpha;x ) \cap C^{s}( \alpha;x ) = \{0\}$,  for any $x \in R^+ \cup R^-$;
\item  $\forall x \in R^+$ we have that $D_{x}T ( C^{u}( \alpha;x )  )\subset C^{u}( \alpha_0;T(x) )$, and 
$$|D_{x}T(v)| \geq \sigma |v|, \mbox{ for all } v \in C^{u}( \alpha;x );$$
\item  $\forall  x \in R^-$ we have that $D_{x}T (  C^{s}( \alpha;x ) ) \subset C^{s}( \alpha_0;T^{-1}(x) )$, and 
$$|D_{x}T^{-1}(v)| \geq \sigma |v|, \mbox{ for all } v \in C^{s}( \alpha;x ).$$
\end{itemize}

It is not difficult to see that if a map $T$ has a hyperbolic cone structure then any  
smooth curve  $\gamma:[0,1] \to R^+$ such that $\gamma'(t) \in C^{u}( \alpha; \gamma(t) )$ for all $t$, then we have that $| T\circ \gamma | \geq  \sigma | \gamma |$, where $|\gamma|$ denotes its length. 

\begin{figure}
\includegraphics*[width=3in]{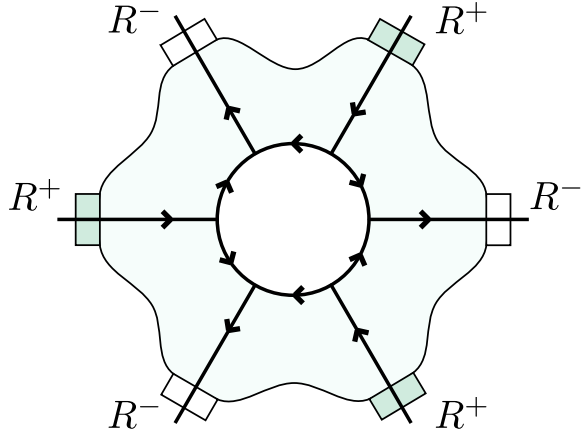}
\caption{Boundary map of degree $3$.}
\label{fig:boundarymap}
\end{figure}

Let us introduce the definition of boundary map, which summarizes the geometric and analytic description of the properties we are interested in. A geometric description of a boundary map is depicted in figure \ref{fig:boundarymap}.
\begin{defi}
Let $r>0$, and let $F:\Cyl{r}\to\Cyl{r}$ be a $C^2-$diffeomorphism. $F$ is a boundary map of degree $p\geq1$ if, under a smooth change of coordinates, the map $F$ satisfies the following properties:
\begin{enumerate}
\item The set of fixed points of $F$ is $\{ ( \pi j /p , 0 )  \, | \, 0 \leq j \leq 2p \}$, and each one is a hyperbolic saddle.
\item The map $F|_{ \SS^1 \times \{0\}}$  is Morse-Smale, and each stable or unstable manifold not contained in $\SS^1 \times \{0\}$ is contained in a line $\ell_{j} = (\pi j / p, \cdot)$, for some $0 \leq j \leq 2p$. 
\item There are real numbers  $0 < r_{-} < r_{+} < r$ and $ \epsilon_{j}>0$, for $0 \leq j \leq 2p$, such that, for each $0 \leq i \leq p$, the following  sets: 
$$R_{i}^{+} =   \left(\frac{2j\pi}{p} - \epsilon _{2j} , \frac{2j\pi}{p} + \epsilon _{2j} \right) \times [r_{-}, r_{+ } ) \smallsetminus \ell_{2j}, $$
$$R_{i}^{-} =   \left ( \frac{(2j+1)\pi}{p}    - \epsilon _{2j+1},  \frac{(2j+1)\pi}{p} + \epsilon _{2j+1}  \right) \times [r_{-}, r_{+ }) \smallsetminus \ell_{2j+1},$$
satisfiy that:  for any $x \in R^+:= \bigcup R_{i}^{+}$ there is a well defined first positive integer $\hat n (x)$ such that $F^{\hat n(x)}(x) \in R^-:= \bigcup R_{i}^{-}$.
\item The map  $T:= F^{\hat n(x)}(x)$ has a hyperbolic cone structure.
\end{enumerate}
\end{defi}

\begin{prop} \label{modelolocal}
Let $p \geq 1$ be an integer number and let $\lambda >1$. If $G \in \diff^2(\Cyl{1}, \SS^1)$ is $C^2-$close enough to the map:
\begin{equation} \label{themap}
F(x,y)  = \left( \frac{2}{p}  \arctan( \lambda^2 \tan(  px/2   ) ) , y \sqrt{ \frac{\cos^2(px/2)}{\lambda^{2} }  + \lambda^2 \sin^2(px/2) } \right),
\end{equation}
then $G$ is a boundary map of degree $p$.
\end{prop}

In order to obtain a proof for Proposition \ref{modelolocal} we need to prove several lemmas. 
Without lost of generality, we will assume $p=2$; similar calculations can be performed for other $p \geq 1$.

Let $\lambda >1$.
The first observation is that the map in (\ref{themap}) is the time-one map of a smooth flow coming from a linear one on the plane after we perform a polar blow up at the origin. 
Let us consider the linear flow $L^t(u,v)= (\lambda^{-t}u,\lambda^t v)$, for $t\in \RR$, and consider the right side of the plane $\RR^2$. 
The horizontal foliation and the vertical foliations are invariant under $DL^t$, the first one contracted and the second expanded.
Moreover, given $\alpha_{0}>0$, there is $0< \rho < 1$ such that, for any $\alpha < \alpha_{0}$ the cone-field:
\begin{equation} \label{linearcones}
C(\alpha; (u,v))^{\perp} = \{w \colon   |\sl^{\perp}(w)| < \alpha \},
\end{equation}
is invariant by $DL^t$, that is,
$DL^t( C(\alpha; (u,v))^{\perp} ) \subset C(\rho^{t} \alpha; L^t(u,v))^{\perp},$
for any $(u,v) \in \RR^2$ and any $t \geq 0$. And the cone-field $C(\alpha; (u,v)) = \{ w \colon  |\sl(w)| < \alpha \}$, satisfies that
$DL^{-t}( C(\alpha; (u,v))  ) \subset C(\rho^{t} \alpha; L^{-t}(u,v)), $ for any $(u,v) \in \RR^2$ and any $t \geq 0$.

Then we perform a polar blow-up of $L^t$ at the origin. This procedure gives us flow in $(-\pi/2,\pi/2)\times (0,r]$ with the following expression:
\[F^t(x,y) = \left(  \arctan(\lambda^{2t} \tan(x) ) ,  y \sqrt{ \lambda^{-2t}\cos^2(x) + \lambda^{2t} \sin^2(x)} \right),\]
that can be extended to a $C^{\infty}-$flow in $[-\pi/2,\pi/2]\times [0,r]$. In order to simplify the notation below, let us rewrite $F^t(x,y) = (h(x,t), y e(x,t))$. 

Both invariant foliations of $L^t$ are transformed into two invariant foliations which are tangent to the vectors: 
\begin{equation} \label{invariantVECTORS}
 \sigma^{u}(x,y) = \left(  \cos(x)  ,  y \sin(x) \right) \text{ and } \sigma^s(x,y) = \left( - \sin(x) , y \cos(x) \right). 
\end{equation}
These two vector fields are obtained by $y \, D_{B^{-1}(x,y)}B (\hat e_{i})$, for $i=2,1$, respectively, where $B(u,v) = (x,y)=(\arctan(v/u), \sqrt{u^2+v^2})$ is the polar change of coordinates.
In the same fashion, the cones $C(\alpha; \cdot)^\perp$ and $C(\alpha; \cdot)$ are transformed into another couple of $DF^t-$invariant cone-fields, say $K^u$ and $K^{s}$, respectively.

It is important to notice that for any fixed $x \neq 0$ we have that: 
\begin{equation} \label{limiteangulo}
\lim_{y= 0} \angle( \sigma^u(x,y), \sigma^s(x,y) ) = 0.
\end{equation}
Hence the angle between the cones $K^u$ and $K^s$ is not bounded from below. 
However, for any $y>0$ the limit in (\ref{limiteangulo}) is $\pi/2$ when $x\to 0$. 

\begin{lem}
Let $U$ be a neighborhood of the origin.
For any $\mu > 0$ there is $\mu_{*}>0$ such that if $y \geq \mu x^2$ and $(x,y) \in U$, we have that
\[\angle ( \sigma^u(x,y), \sigma^s(x,y) ) \geq \mu_{*}. \]
\end{lem}

\begin{proof}
Let us consider new coordinates $(w,s)$ in the plane given by $\tilde B(w,s)=(w,sw)$ and whose inverse is $\tilde B^{-1}(x,y)= (x, y/x)$. In this coordinates the flow $F^t$ induces a new flow given by $H^t(w,s) = \tilde B^{-1} \circ F^t \circ\tilde B(w,s)$, for $t\in \RR$. This flow has the following expression:
\[ H^t(w,s) = \left(  h(w,t), s w \frac{ e(w,t) }{ h(w,t) }    \right). \]

The flow $H^t$ can be extended to the line $[w=0]$, by l'H\^opital rule. In fact, the limit
\[ \lim_{w=0} \frac{w}{h(w,t)} = \left(  \left. \frac{\partial h(w,t) }{ \partial w}\right|_{w=0} \right)^{-1}  = \lambda^{-2t} \] 
is well defined since $h(w,t)$ is $C^{\infty}$ in $w$. Of course, it is needed only that $h$ is $C^2$, to obtain a resulting flow of class $C^1$.
This procedure do not make any change in the dynamics of the linear flow $A^t$ first considered; we have only added one more line $[w=0]$ to the space.

Again, there are two invariant foliations of $H^t$, tangent to the vector fields:
\[\tilde \sigma ^u (w,s) = (w \cos(w), -s \cos(w) + s w \sin(w)), \]
\[\tilde \sigma ^s (w,s) = (- w \sin(w),  s \sin(w) + s w \cos(w)), \]
which come from $\sigma^s$ and $\sigma^{u}$, as in (\ref{invariantVECTORS}).  
The angle between these two cone fields, in a region where $w$ is close to zero, has the following expression:
\begin{equation} \label{angulodos}
\angle( \tilde \sigma^u(w,s), \tilde \sigma^s(w,s) )  =  \frac{w}{2s} + O(w^3).
\end{equation}
Hence, if $s \geq  \mu w$ for some $\mu >0$ then the angle in (\ref{angulodos}) is bounded from below by $(2\mu)^{-1}$. That is, the angle is bounded from below for points above a line of slope $\mu$ through the origin, and tends to zero below such a line.
Notice that, in coordinates $(x,y)$, condition that $s \geq  \mu w$ is that  $y \geq \mu x^2$. Hence, there is some $\mu_{*}>0$ such that satisfies the Lemma.   
\end{proof}

Observe that we can translate the cone-fields $K^{u,s}$ into two cone-fields in the $(w,s)$-plane. These cone-fields are preserved by any $C^1-$perturbation of $H$.
Moreover, any $C^1-$perturbation of $H$ comes from a $C^2-$perturbation of $F$, and viceversa. So, mixing this property with the fact that above some fixed $y>0$ the angle between  $K^{u}$ and $K^s$ is bounded from below, we obtain the following Corollary: 
\begin{cor}
Let $\mu>0$, if $G \in \diff^2(\Cyl{1}, \SS^1)$ is a map $C^2-$close enough to the map $F$, then if  $x$ is small and $y \geq \mu_{*} x^2$,  then
$$D_{(x,y)}G( K^{u}(\cdot ;(x,y) ) ) \subset K^{u}(\cdot ;G(x,y) ), $$
and the analogous property for the stable cone-field $K^s$.
\end{cor}

The map $F$ can be $C^\infty-$linearized at the origin, even that the eigenvalues of $D_{0}F$ are ressonant. This linearization is given by 
$ \Psi (x,y) = (\tan(x), y \cos(x) )$  and $\Phi(u,v)= (\arctan(u), v \sqrt{1+u^2})$. Denote by 
\[L_{0} (u,v)= D_{0}F(u,v) =  (\lambda^2u, \lambda^{-1}v ).\]
If $G$ is $C^1-$close enough to $F$, then there is a neighborhood $U_{G}$ of the origin where $G$ is $C^1-$linearizable.
In fact, Hartman proved a Theorem about $C^1-$linearization of hyperbolic saddles, in any dimension, with some geometric restriction on the eigenvalues (see \cite[pag. 235]{Hartman}).   
Although the existence of a $C^1-$linearizing neighborhood, in general, depends on the nature of the resonance of the eigenvalues (see \cite{Sell}, for instance), the case of the plane is special and Hartman's result guarantees that any hyperbolic saddle in the plane is $C^1-$linearizable. 

Consider the partition of the interval $(0,\pi/2)$ by intervals of the form $(x_{n},x_{n+1}]$, induced by the sequence $x_{n} = h (\arctan(\lambda^{-1}),n)$, for $n \in \ZZ$. Notice that $\tan(x_{n}) =  \lambda^{2n-1}$, for any $n$.
The derivative of $F$ at any point $z=(x,y)$ has the form:
\[ D_{z}F=\begin{pmatrix}  \left. \frac{ \partial  }{\partial x} F_{1}\right |_{z}  & \left. \frac{ \partial  }{\partial y} F_{1}\right |_{z}  \\ 
\left. \frac{ \partial  }{\partial x} F_{2}\right |_{z}  & \left. \frac{ \partial  }{\partial y} F_{2}\right |_{z}  \\  \end{pmatrix} =  \begin{pmatrix}  h'(x) & 0 \\ ye'(x) & e(x) \\  \end{pmatrix}, \]
where $h(x)=h(x,1)$ and $e(x)=e(x,1)$. So, we can state the following Lemma.

\begin{lem}  \label{boundsderivative}
If $G=(G_{1},G_{2}):M \to M$ is a $C^2$ map which is $C^2-$close to $F$, 
then there is number $\tilde \lambda >1$ such that:
\begin{enumerate}
\item $1 < \left. \frac{ \partial  }{\partial x} G_{1}\right |_{z}  <\tilde \lambda^2$, for $z \in (0,\tilde x_{0}) \times \RR^+$;
\item $1/\tilde  \lambda^2 < \left. \frac{ \partial  }{\partial x} G_{1}\right |_{z} <1$, for $z \in (\tilde x_{0}, \pi/2)\times \RR^+$;
\item $1/\tilde \lambda < \left. \frac{ \partial  }{\partial y} G_{2}\right |_{z}    < \tilde \lambda$, for $z \in [0,\pi/2] \times \RR^+$;
\item The function $\left. \frac{ \partial^2  }{\partial y \partial x} G_{2} \right |_{z}$ has a unique critical point in $[0,\pi/2]$ at some point close to $x_{0}$, and this point is a maximum.
\end{enumerate}
\end{lem}

\begin{proof}
Denote by $x_{0} = \arctan(1/ \lambda ) \in (0,\pi/2)$. An easy computation gives all inequalities for $F$ itself, that is:
(1) $1 < h'(x) <\lambda^2$ for $x \in (0,x_{0})$;
(2) $1/\lambda^2 < h'(x) <1$ for $x \in (x_{0}, \pi/2)$;
(3) $1/\lambda < e(x) < \lambda$ for $x \in [0,\pi/2]$;
and (4) the point $x_{0}$ is the unique critical point of $e'(x)$ and is a global maximum.  
Since $G$, in particular,  is $C^1-$close to $F$, it is easy to see that inequalities (1)-(3) hold for $G$ too. Notice that (4) is not true for $C^1-$perturbations. However, since $G$ to be $C^2-$close to $F$ we have that, for $z = (x,y)$, 
\[ \left |   \left. \frac{ \partial^2 G_{2} }{\partial y \partial x}  \right |_{z} - \left. \frac{ \partial^2  F_{2}}{\partial y \partial x}   \right |_{z}   \right |
= \left |   \left. \frac{ \partial^2 G_{2} }{\partial y \partial x}  \right |_{z} - e'(x) \right |\]
is uniformly small, so the unique maximum of $e'(x)$ is preserved.
\end{proof}

\begin{lem} \label{lemakickup}
Let $G$ be map which is $C^2-$close to $F$, and let $U=U_{G}$ be a $C^1-$linear\-iza\-tion neighborhood of $G$. Then there are $N\leq -1$ and $Y>0$, such that the rectangle $R =[0,x_{N}] \times [0,Y] \subset U_{G}$, and there is $\tilde \lambda =\min\{|\!| D_{0}G|\!| , \lambda  \}>1$, such that:
\begin{enumerate}
\item If $z \in R \smallsetminus F(R)$, $F^j(z) \in R $ for all $j \in \{1, \ldots, m\}$, and $F^{m+1}(z) \notin R$, and if $m$ is large enough then
\begin{equation*} 
K_{m}:= D_{G^m(z)} \Phi [ D_{z}G^m( K^u(\cdot; z) ^+ )]
\end{equation*}
is contained in the positive quadrant. Moreover, there is a positive constant $c$, close to 1, such that for any $v \in K_{m}$ we have that 
\[\sl(v) > \eta(m):= c Y \tilde \lambda^{2N - m - 4} - \beta \tilde  \lambda^{-3m} >0.\]
\item If $\tilde z=(\tilde x, \tilde y) \in (x_{N-1},x_{N}]\times[0,1]$, we have that $F^{|N| +1 }(\tilde z) \in [x_{0},x_{1}] \times [0,1]$, and $\sl(v) >0$, for every $v \in D_{\tilde z}G^{|N| +1 }(K_{m})$.
\end{enumerate}
\end{lem}

\begin{proof}
Let $U$ be a small open neighborhood of the origin, contained in the neighborhood of $C^1-$linearization of $F$, 
and consider $Y>0$ and $N\leq -1$ such that the rectangle
$R= [0,x_{N}] \times [0,Y]$ is contained in $U$. 
There is $\beta >0$ such that, if $z \in R \smallsetminus F(R)$ then 
\[D_{z}\Phi(K^{u} (\cdot ; z) ) \subset C(\beta; \Phi(z) ),\]
where $C(\beta; \Phi(z) ) = \{v \, | \, |\sl(v)| < \beta \}$.
On the other hand, the map $L_{0}$ is linear, so for any $(u,v)$ and $n\geq 0$ we have that:
\[ L_{0}^n(C(\beta; (u,v) )) = C(\lambda^{-3n}\beta ; L_{0}^n(u,v) ).\]
Hence, if $F^j(z) \in R$, for all $j \in \{1, \ldots, m\}$ and $m\geq 1$, we have that: 
\begin{equation} \label{spreadcontraction}
D_{F^m(z)} \Phi [ D_{z}F^m( K^u(\cdot; z) )] \subset C(\lambda^{-3m}\beta ; \Phi (F^m(z) )  ).
\end{equation}

Take a point $z = (x,y) \in R \smallsetminus F(R)$ that $F^j(z) \in R$, for $j \in \{0,\ldots m\}$ and assume $F^{m+1}(z) \notin R$. Denote by $(\tilde x,\tilde y) = F^m(z)$. This implies that $x_{N-1} \leq \tilde x \leq x_{N}$, and hence $\tan(\tilde x) \in [\lambda^{2N-3}, \lambda^{2N-1}]$. On the other hand, there is a constant $c>0$, close to 1, such that $y \geq c \lambda^{-1}Y$, since the unstable foliation of $F$ is $C^1-$close to the horizontal in this neighborhood. Therefore, $\tilde y \geq c \lambda^{-m-1}Y$. If we write $w=\Phi(\tilde x, \tilde y)$, then $K_{m}$ is a cone around the vector $\check{e} = D_{w}\Psi(\hat e_{1}) = \cos(\tilde x)( \cos(\tilde x), \sin(\tilde x) )$, and the slope $\sl(\check{e}) = \tilde y \tan(\tilde x) \geq c Y \lambda^{-m-1} \lambda^{2N-3}  >0$.

Now take any $v \in K_{m} \subset D_{F^m(z)} \Psi [ C( \beta  \lambda^{-3m}; \Phi(F^m(z) )  )^+  ]$, then (\ref{spreadcontraction}) imply that
\[\sl(v) \geq \sl(\check{e}) - \beta \lambda^{-3m} = cY \lambda^{2N - m - 4} - \beta \lambda^{-3m}.\]
So, if $m \geq 1$ is large enough, we obtain item (1) for $F$. 
Let $G$ be a map which is $C^2-$close to $F$, and let $U_{G}$ be a $C^1-$linearization neighborhood. 
If $N\leq -1$ and $Y$ satisfy that $R \subset U_{G}$, the previous estimates gives a proof of item (1) for such a $G$.

In order to get item (2), we have to guarantee that the iterates of the cone $K_{m}$ remain on the positive quadrant, in the following $|N|$ iterations, until the point reaches the subset $(x_{-1},x_{0}] \times [0,1]$.
For this, bounds for partial derivatives of $F$, up to second order, given in Lemma \ref{boundsderivative}, are the key. 

Recall that the point $(\tilde x, \tilde y) =F^m(z)$ and take any vector $(v_{1},v_{2}) \in K_{m}$. 
So, $F(\tilde{x},\tilde y) \in (x_{N},x_{N+1}]\times [0,1]$. If we let 
\[(w_{1},w_{2}) =D_{(\tilde x, \tilde y)}F(v_{1},v_{2}) =  (h'(x) v_{1} , y e'(x) v_{1} + e(x) v_{2} ),\]
we can bound  $\sl(w_{1}, w_{2})$, away from zero. For this we use that 
$\sl(v_{1},v_{2}) > 0$, for any $v \in K_{m}$, and $\min\{e(x)\} \geq \lambda^{-1}$. Moreover, since $x_{N-1} \leq \tilde x$, then
\[0< \tilde y e'(x_{N-1}) \leq \tilde y e'(\tilde x ).\] 
Hence, we have:
\[\sl (w_{1},w_{2})  =   \frac{\tilde y e'( \tilde x) v_{1} + e( \tilde x) v_{2}}{h'( \tilde x) v_{1}}  
 \geq  \frac{\min\{ \tilde y  e'(\tilde x)   + e(\tilde x) ( v_{2}/v_{1} ) \}}{\max\{   h'(\tilde x)      \}} \geq \frac{ \sl(v_{1},v_{2})}{\lambda^{3}} >0. \] 
We can repeat this computation $|N|+1$ times, until we reach the last interval $(x_{0},x_{1}]$, proving that
\[\sl ( D_{(\tilde x, \tilde y)} F^{|N|+1} (v) ) >   {\lambda^{ -3(|N|+1)}} \sl(v )> 0,\]
for any vector $v \in K_{m}$. Notice that $x_{N-1}\leq x_{1}$, also. 
This gives the proof of item (2) for $F$. Finally, Lemma \ref{boundsderivative}, allows us to extend this bounds
for the perturbation $G$.
In fact we obtain that 
\[\sl ( D_{(\tilde x, \tilde y)} G^{|N|+1} (v) ) > {\tilde \lambda^{- 3(|N|+1)}}  (cY \tilde \lambda^{2N-m-4}  - \beta \tilde \lambda^{-3m} ) > 0,\]
for any $v \in K_{m}$.
\end{proof}

Now we can give a proof of Proposition \ref{modelolocal}.
\begin{proof}[Proof of Proposition \ref{modelolocal}]
Let $\lambda >1$ and $p\geq 1$. 
The map $F: \C_{\infty} \to \C_{\infty}$ defined in (\ref{themap}) is a $C^{\infty}-$diffeo\-mor\-phism with $2p$ hyperbolic fixed points on the boundary. 
These periodic points have their stable and unstable manifolds contained in $\SS^1 \times \{0\}$ coincide alternately. 
We shall restrict to the region bounded by $P_{\pm} = (\pm \pi/p,0)$, that is $[P_{-}, P_{+}]\times \RR$. 
Inside this region there is another hyperbolic fixed point $P_{0} = (0,0)$. 
We can compute easily the derivative of $F$ at these points:
\[D_{P_{0}} F = (D_{P_{\pm}} F)^{-1} = \begin{pmatrix} \lambda^{2} & 0 \\ 0 & \lambda^{-1} \\ \end{pmatrix}.\]
In particular, $F$ is Morse-Smale restricted to the boundary of $\C_{\infty}$.
Since the map $F$ is symmetric with respect to the line $[x=0]$, we shall focus only on $[0,P_{+}] \times [0,1]$, see figure \ref{fig:passage}.

\begin{figure}
\includegraphics*[width=4in]{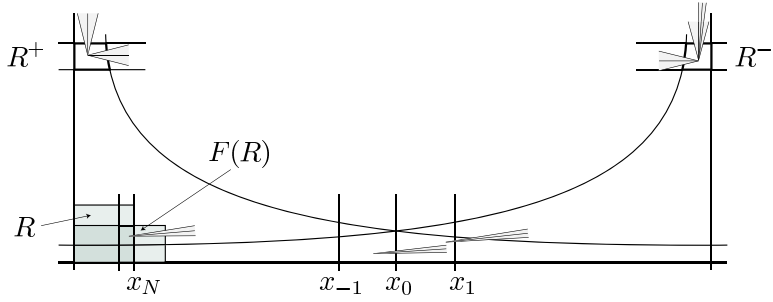}
\caption{Transitions through the boundary}
\label{fig:passage}
\end{figure}

In the following we will consider the case $p=2$ which simplifies the notation; a similar argument can be performed by other $p\geq 1$. 
Consider $G \in \diff^2(\Cyl{1}, \SS^1)$ which is $C^2-$close to the map $F$, with $p=2$, that satisfies all previous Lemmas.
Let $\eps_{+}$, $\eps_{-}$ and $r_{+}$ be three small positive numbers and consider two regions:
\begin{equation*}  
R^{+} = R^{+}(\eps_{+}, r_{+}) = [0,\eps_{+}] \times[\lambda^{-1} r_{+},  r_{+}) 
\end{equation*}
$$R^{-} = R^{-}(\eps_{-}, r_{+}) =  [\pi/2 - \eps_{-}, \pi/2] \times[\lambda^{-1} r_{+},  r_{+}),$$
In order to prove that the map $G$ is a boundary map of degree $2$, we need to choose these numbers adequately such that 
the transition map $T:R^+\to R^-$, defined as $T(x,y) := F^{\hat{n}}(x,y)$, for $\hat{n} = \hat{n}(x,y) \geq 1$ that $G^{\hat{n}} (x,y) \in R^{-}$, has a hyperbolic cone structure. 
Observe that if $\eps_{+}$ and $\eps_{-}$ are chosen small enough, then $\hat n$ is well defined and unique.
Also it is true that $\min\{ \hat n(x,y) | (x,y) \in R^+\} \to \infty$, uniformly on $(x,y)$, as $\eps_{+}$ tends to $0$. 
 
By construction of $F$, the cone-fields defined in (\ref{linearcones}), for some fixed $\alpha>0$ are mapped to a couple of cone-fields $K^u(\alpha, z)$ and $K^s(\alpha, z)$, restricted to points in $R^{+} \cup R^-$, and they are almost orthogonal. In the following we shall prove that these cones provide of a hyperbolic cone structure for $T$.
 
Let $U_{G}$ be the $C^1-$linearization neighborhood of $G$. 
Let $N\leq -1$, $Y>0$ and $R$ given in Lemma \ref{lemakickup}. 
Take any point $z \in R^+$.
If $\eps_{+}$ is small enough, then we have that there is an integer $n_{1}>0$ such that $F^{n_{1}}(z) \in R \smallsetminus F(R)$. 
Perhaps considering a smaller $\eps_{+}$, we can assume also that $ F^{j+n_{1}}  (z) \in R$ for $j \in \{0, \ldots, m\}$ and $F^{j+1+n_{1}}  (z)$, for $m$ large enough such that item (1) of Lemma \ref{lemakickup} holds.
Hence we have that for any vector $v \in D_{z}G^{j+1+n_{1}}(K^u(\alpha;z))^+$ satisfies that: 
\[\sl(v) > \eta(m) >0.\] 

In order to guarantee that this last cone lands inside some $K^u( \cdot ; G^{\hat n}(z) )$, at $R^-$, we need to verify it for $F$ first.
Observe that if we let $\eps_{+} > 0$ to be small enough then, for a fixed small $\eps_{-}$, we have that for any $z \in R^{+}$ the value~$n_{2} := \hat n(x,y) - ( m + |N| ) > m^*$, for any given arbitrarily large $m^* > 0$.

Let $(x,y) =  F^{j+n_{1}}(z)$. The derivative of $D_{(x,y)}F^{n}$, for any $n \in \ZZ$ is
 \begin{equation*}  
  D_{( x, y)}F^{n} =    \begin{pmatrix}  
       h_{n}'(  x) & 0 \\
       \tilde y e_{n}'(  x) & e_{n}(  x) \\
    \end{pmatrix}
 \end{equation*} 
where $h_{n}(x) = h(x,n )$ and $e_{n }(x) = e(x,n)$.
Given a vector $(v_{1},v_{2})$ that $\sl^{\perp}(v_{1},v_{2}) \geq \eta(m)^{-1} \sim \lambda^m$ then:
\begin{eqnarray*}
\sl^{\perp}( D_{(x,  y)}F^{n_{2}}( v_{1},v_{2} ) ) = \frac{h'_{n_{2}}( x) v_{1}} {  y e'_{n_{2}}( x) v_{1} + e_{n_{2}}( x) v_{2}  }    \\
\leq  \frac{h'_{n_{2}}( x)  } {   e_{n_{2}}( x)    }\sl^{\perp}(v_{1}, v_{2}) 
\leq \frac{h'_{n_{2}}( x)  } {   e_{n_{2}}(x)    } \frac{1}{\eta(m)} . 
\end{eqnarray*}
Since $y e'_{m}( x)v_{1}>0$. Moreover, $ h'_{m}(x) / e_{m}(x) \sim \lambda^{-2m} $, hence, there is $m^*>0$ such that 
$$\sl^{\perp}( D_{( x,  y)}F^{n_{2}}( v_{1},v_{2} ) ) < \alpha_{0} < \alpha.$$
Therefore, we have for any $z\in R^+$ that:
$$D_{z}T ( K^u(\alpha, z) ) = D_{(x,y)}F^{\hat n} (  K^u(\alpha, z)  ) \subset C^{\perp}(\alpha_{0}, F^{\hat n}z) \subset K^u(\alpha, F^{\hat n}z).$$
It is not difficult to see that vectors inside the cone $K^u (\alpha, \cdot)$ are expanded by the derivative of $T$, for $F$.
Finally, Lemma \ref{boundsderivative} allows us to extend these computations to the derivative of $G$. 
On the other hand, notice that the cone-field $K^s(\alpha, \cdot)$, on points $R^-$ is preserved (and vectors inside are expanded) by $DT^{-1}$. In fact, the same computations for $0<\lambda <1$ instead of $\lambda >1$, yield to this result. This finishes the proof of Proposition \ref{modelolocal}.
\end{proof}
 
A straightforward consequence of Lemma \ref{blowupconstruction} and Proposition \ref{modelolocal} is the following.
\begin{prop}\label{keyprop}
If $f$ is a blow-up of a pseudo-Anosov map, there is $\mathcal{N}_{2}\ni f$, a $C^2-$neighborhood, such that any $g\in \mathcal{N}$, restricted to a neighborhood of a component of the boundary is a boundary map. 
\end{prop}

It is important to mention that as a consequence of the proof of Proposition \ref{modelolocal}, a map $G:\C_{r} \to \C_{r}$ which is $C^2$-close enough to $F$ induces a map $T:R^+ \to R^-$ which has a kind of \emph{Markov property}: for any curve $\gamma$ crossing $R^+$ along the horizontal direction and whose tangent vectors are contained in the cone $K^{u} (\alpha)$ is that $T(\gamma)$ crosses $R^-$ along the vertical direction. Moreover, the function $\hat n : R^+ \to \NN$ induces a partition $\mathcal{P}^+ = \{ \Delta_{t}^+:= \hat n ^{-1}(t) | t \in \NN \}$  of $R^+$ where the map $T:R^+ \to R^-$ is continuous, and each $\Delta_{t}^+$ is foliated by stable and unstable leaves which are actually contracted and expanded by $T$, and more, $T(\Delta_{t}^+)$ crosses $R^-$ along the vertical direction.

\subsection{Transitivity}
Let $(S, \partial S)$ be a surface with boundary of genus $g \geq 1$. Although the following construction can be performed for any number of connected components of $\partial S$, notation becomes simpler if we assume that $\partial S $ has only one connected component, and we do so. Let $S_{0}$ be the  surface without boundary obtained after collapsing $\partial S$ to a point $\sigma$. The configuration $(\sigma, p)$ with $p = 4(g-1)+2 \geq 2$ satisfies (\ref{mscondition}) of Theorem \ref{ms}, so there is a pseudo Anosov map $f_{0}:S_{0} \to S_{0}$ with dilatation $\lambda >1$ and one $p$-prong singular point $\sigma$.  Moreover, we can assume that there is a hyperbolic fixed point $w \neq \sigma$, in $S_{0}$, whose stable and unstable manifolds are dense. This is true, perhaps considering some iterate of $f_{0}$ instead, since there are plenty of hyperbolic periodic points. 

Let $f \in \diff^{\infty}(S,\partial S)$ be the blow-up of the pseudo-Anosov at $\{\sigma\}$. We already know that this map is transitive, since we have not changed the orbit structure of the original pseudo-Anosov $f_{0}$. Now we can give a proof of Theorem \ref{p-a-thm}.

\begin{proof}[Proof of Theorem \ref{p-a-thm}]
Let $W$ be a small neighborhood of $\partial S$ on $S$, and let $g \in \mathcal{N}_{1} \cap \mathcal{N}_{2}$, from  Lemma \ref{unifhyp} and Proposition \ref{keyprop}.
The neighborhood $W$ is the union of $p$ sectors $\Gamma_{i}$ corresponding to the sectors between unstable prongs of $f_{0}$ landing at $\sigma$.
Denote by $U =S\smallsetminus W$. Observe that Lemma \ref{unifhyp} guarantees that $\Lambda = \bigcap_{n \in \ZZ}g^n(U)$ is a hyperbolic maximal invariant set of $g$.  
In fact, $g|_{\Lambda_{g}}$ is transitive and $\Lambda_{g}= H_{q }$, where $q  \in U$ is a hyperbolic fixed point of $g$.
Proposition \ref{keyprop} states that 
on each sector $\Gamma_{i}$ there are input and output sections $R^+_{i}(\eps_{i+},r_{i})$ and $R^-_{i}(\eps_{i-}, r_{i})$ which satisfy the properties of a boundary map for $T: R^+ \to R^-$ induced by $g$, where $R^+:=\cup R^+_{i}$ and $R^-:=\cup R^-_{i}$. 
The input section $R^+$ is a finite union of rectangles. Without loose of generality, we can assume that the horizontal boundary of each $R^+_{i}$ is given by pieces of the unstable manifold of $q $, and one of their vertical boundaries is given by the stable manifold of some a fixed point at the boundary.  Analogously, the output section $R^-$ is a finite union of rectangles whose horizontal boundaries are given by pieces of the stable manifold of $q $ and one of its vertical boundaries is given by the local unstable manifold of a fixed point of the boundary.

Now, let $A$ be any non-empty open set of $S$, then we claim that the sets 
\[ \{ j \geq 0 \, | \, g^j(A) \cap R^+ \neq \emptyset \} \text{ and } \{ j \geq 0 \, | \, g^{-j}(A) \cap R^- \neq \emptyset \}\]
are unbounded.
To proof this claim for the former set, let us assume, by contradiction, that $g^j(A) \subset U$ for all $j \geq 0$. This implies that $A \subset W^s(\Lambda_{g})$. Hence, the Shadowing Lemma and the fact that $f$ is transitive in $\Lambda_{g}$ imply that there is a point $z_{1} \in A$ and a sequence $n_{j} \to \infty$ such that $g^{n_{j}}(z_{1}) \to q$. 
Moreover, $W^s(q) \cap A \neq \emptyset$. In fact, consider a small curve $\gamma$ through $z_{1}$ inside the unstable cone field $K^{u}(\alpha)$; expansion of vectors inside these cones guarantee that $g^{n_{j}}(\gamma)$ crosses $W_{\text{loc}}^s(q)$, for some $n_{j}$. 
Now observe that Palis' Inclination Lemma guarantees that there are $\gamma_{0} \subset \gamma$  and $m \geq 0$ such that $g^m(\gamma_{0}) \cap  R^+ \neq \emptyset$. 
Recall that horizontal boundaries of $R^+$ are contained in compact pieces of the unstable manifold of $q$.
Hence, there are points in $A$ that escape from $U$, since $ R^+ \subset S \smallsetminus U$. This contradiction proves the claim. The same argument implies that the latter set is unbounded.

Now we can prove that for any two not empty open sets $A$ and $B$  of $S$, there is $n\in \ZZ$ that $g^n(A) \cap B \neq \emptyset$.
Without lost of generality, we can assume that both open sets are contained in $U$, since no open set remains all its iterates inside $W$.
As a consequence of the previous claim, there is a point $a \in A$   that the forward  orbit intersects $R^+$ infinitely many times, so there is some $n_{1}\geq 0$ and there is a curve $\gamma_{u} \subset A$ inside the cone field $K^{u}(\alpha)$, such that $g^{n_{1}}(\gamma_{u})$ crosses $R^+$ along the horizontal direction. 
On the other hand, there is a point $b \in B$ such that its backward orbit intersects $R^-$ infinitely many times. 
Hence, there is an integer $n_{2} \geq 0$ and a curve $\gamma_{s}  \subset B$, inside the cone field $K^s(\alpha)$, such that $g^{- n_{2}}(\gamma_{s})$ crosses $R^-$ along the horizontal direction. Finally, $T(g^{n_{1}}(\gamma_{u}))$ is a vertical curve in $R^-$ and hence $T(g^{n_{1}}(\gamma_{u})) \cap g^{-n_{2}}(\gamma_{s}) \neq \emptyset$, and we are done.
\end{proof}

Observe that the previous claim finished the proof of Theorem \ref{p-a-thm} and therefore Theorem \ref{cedos}. To conclude Theorem  \ref{cedos-hc}, observe that the proof of Theorem \ref{p-a-thm} implies that any forward (backward) iterate of an open set intersects transversely the stable (unstable) manifold of the periodic point with orbit far from the boundary, and this implies that there are homoclinic points inside any open set. 

To conclude Theorem \ref{nonhyp} observe that the presence of boundary immediately prevent the map to be globally hyperbolic and they can not have homoclinic tangencies since otherwise, its unfolding would create (generically)  either sinks or repelling and therefore, the maps would not be robustly transitive.

\bibliographystyle{plain}

\end{document}